\documentclass[12pt]{article}
\usepackage{amsmath}
\usepackage{amssymb}
\usepackage{amsthm}
\usepackage{amsfonts}
\usepackage{amscd}
\usepackage[mathscr]{eucal}
\newcommand{\Z} {{\mathbb  Z}}
\newcommand{\Q}{{\mathbb  Q}}

\begin{document}
\parindent  25pt
\baselineskip  10mm
\textwidth  15cm    \textheight  23cm \evensidemargin -0.06cm
\oddsidemargin -0.01cm

\title{ { Circulant Digraphs Integral over Number Fields
 }}
\author{\mbox{}
{Fei Li}
\thanks{ \quad E-mail:
cczxlf@163.com } \\
(School of Statistics and Applied Mathematics, \\ Anhui University of Finance and Economics,\\
 Bengbu City, 233030, Anhui Province, P.R.China)  }

\date{}
\maketitle
\parindent  24pt
\baselineskip  10mm
\parskip  0pt

\par   \vskip 0.6cm

{\bf Abstract} \quad A number field $K$ is a finite extension
of rational number field $ \Q. $  A circulant
digraph integral over $K$ means that all its eigenvalues are algebraic integers of $ K. $ In this paper
we give the sufficient and necessary condition for circulant digraphs which are integral
over a number field $ K. $ And we solve the Conjecture3.3 in [XM] and find it is affirmative.

\par     \vskip  1.3cm

\hspace{-0.6cm}{\bf 1 \ \ Introduction and Main Results}

\par \vskip 0.6cm

Let $ G $ be a finite group and $ S $ be a subset of $ G. $ The Cayley digraph
$ D = D(G, S) $ of $ G $ with respect to $ S $ is a directed graph with vertex
set $ G. $ For $ g_{1}, g_{2} \in G, $ there is an arc from $ g_{1} $ to $ g_{2} $
if and only if $ g_{2}g_{1}^{-1} \in S. $ If $ G $ is the cyclic group $ \Z/n\Z $
with order $ n, $ the Cayley digraph $ D = D(\Z/n\Z, S) = D(n, S) $ is called a
circulant digraph. Its spectrum consists of the eigenvalues of the adjacent matrix of
$ D, $ namely, the roots of the characteristic polynomial
of $ D, $ that is $ spec(D(n, S))=\{\lambda_{0}, \lambda_{1}, \cdots, \lambda_{n-1}\}, $ where
$ \lambda_{r}=\sum _{s \in S}\zeta_{n}^{rs}(1 \leq r \leq n-1) $ with
$\zeta_{n}=e^{\frac{2\pi i}{n}}, $ the primitive $n-$th root of unity(see [Bi], [GR]).
It is obvious that $ \lambda_{r} (1 \leq r \leq n-1) $ are algebraic integers
in the cyclotomic field $ \Q(\zeta_{n}). $ A circulant digraph integral over number field
$K$ means that all its eigenvalues are algebraic integers of $ K. $ The reader is referred
to [F], [N] and [W] for terms left undefined and for the general theory of number fields.

Our main purpose in this paper is to give the sufficient and necessary condition
for circulant digraphs which are integral over a number field $ K. $
And we solve the Conjecture3.3 in [XM] and find it is affirmative. There are a lot of
literature studying integral graphs(see e.g., [Ba][BC][EH][HS][LL][S][Z]).

The main results in the present paper are stated as follows.

 {\bf A.  (see Theorem 1.) }  The circulant directed graph $ D(n, S) $ is integral
over $ K $ if and only if $ S $ is a union of some $ M_{i}'$s.

{\bf B. (see Corollary 1.) }  There are at most $ 2^{r(n, K)} $ basic
circulant digraphs integral over $ K $ on $ n $ vertices.

Using the same notations of
[XM], we have the following corollary.

{\bf C. (see Corollary 2.) } (Conjecture 3.3 in [XM]) Let $ D(\Z_{n}, S) $ be a
circulant digraph on $ n = 2^{t}k(t >2) $ vertices, then $ D(\Z_{n}, S) $
is Gauss integral if and only if $ S $ is a union of the $M'$s.

\par \vskip 1 cm

\hspace{-0.6cm}{\bf 2 \ \ Constructing a partition of \ $ \{1, 2, \cdots, n-1\} $}

\par \vskip 0.4 cm

Let $ L/M $ be a number field extension. We denote $ Gal(L/M) $ the Galois group of L
over M and $ [L:M] $ the dimension of L as vector space over M.
Let $ K $ be a number field, i.e., a finite extension of rational number field $ \Q, $
and $ d(n) $ be the set of all the positive proper divisors of $ n $ and for
$ p\in d(n), $ we write $ n = p g. $ Denote $ F = K \bigcap \Q(\zeta_{n}) $ and for
$ p \in d(n) $ define $ G_{n}(p) = \{g | 1 \leq g < n, gcd(g, n) = p \}. $ It is
clear that $ G_{n}(p) = p G_{g}(1). $ Notice that $F$ is an abelian field and
$ D(n, S) $ is integral over $K$ if and only if $ D(n, S) $ is integral over $ F. $

Let $ H^{\prime} = Gal(\Q(\zeta_{n})/F), H = Gal(F \cdot \Q(\zeta_{g})/F), \
G^{\prime} = Gal(\Q(\zeta_{g})/F \bigcap \Q(\zeta_{g})), $ and $G = Gal(\Q(\zeta_{g})/\Q). $
By Galois theory( see [L], P.266), $ H \cong G^{\prime} \subseteq G. $ Denote $\pi$
the translation operation of $G^{\prime}$ on $ G $ and $ f: G \longrightarrow
(\Z/g\Z)^{\ast} $ the map defined by $ f: \sigma \mapsto a $ with $ \sigma(\zeta_{g})
= \zeta_{g}^{a}. $ The map $ f $ is an isomorphism of groups, by which
$\pi$ also can be considered as an operation on $(\Z/g\Z)^{\ast}. $ For
$ G_{n}(p) = p (\Z/g\Z)^{\ast}, $ we can define an operation $\pi^{\prime}$ of
$G^{\prime}$ on $ G_{n}(p) $ by $ \pi_{g}^{\prime}(x) = p\pi_{g}(x^{\prime}) $ with
$ x = px^{\prime} \in G_{n}(p), g \in G^{\prime} $ and $ x^{\prime} \in (\Z/g\Z)^{\ast}. $
It is easy to see the number $r_{p}$ of all orbits under $ \pi^{\prime} $ is equal
to $[F \bigcap \Q(\zeta_{g}) : \Q]$ and each orbit $M_{i}(1 \leq i \leq r_{p})$ has
$[\Q(\zeta_{g}) :F \bigcap \Q(\zeta_{g})]$ elements. By Galois theory, $ H^{\prime}, $ when
restricted to $ F \cdot \Q(\zeta_{g}), $ is equal to $H$ and $ H \cong G^{\prime}. $ Hence
$\pi^{\prime}$ can be regarded as an operation of $ H^{\prime} $ on $ G_{n}(p), $ under which
there are the same $r_{p}$ orbits. So $ G_{n}(p) $ is the disjoint union of the distinct
$ M_{i}'$s and we can write $ G_{n}(p) = \bigsqcup_{i=1}^{r_{p}}M_{i}. $  It is clear that
$\{1, 2, \cdots, n-1\}$ has a partition written as $\{1, 2, \cdots, n-1\} =
\bigsqcup_{p \in d(n)}G_{n}(p). $ So $\{1, 2, \cdots, n-1\}$ has a new partition
$ \bigsqcup_{i=1}^{r(n, K)}M_{i}, $ where $ M_{i} $ is an orbit for some $ p \in d(n) $
and $ r(n, K) = \sum _{p \in d(n)}r_{p}. $

\par \vskip 0.4 cm

\hspace{-0.6cm}{\bf 3 \ \ Proofs of Main Results}

\par \vskip 0.4 cm

If $ S $ is a union of some $ M_{i}'$s, we have $ \sigma S = S $ for each
$ \sigma \in H^{\prime}$ by the construction of $ M_{i}'$s.
So $\sigma(\lambda_{r}) = \lambda_{r}(0 \leq r \leq n-1), $ which shows that
$ \lambda_{r} \in F \subseteq K $( see [L], P.262) and $ D(n, S) $ is integral over $ K. $ We have
the following theorem.
\par \vskip 0.4 cm

 {\bf Theorem 1. } \ The circulant directed graph $ D(n, S) $ is integral
over $ K $ if and only if $ S $ is a union of some $ M_{i}'$s.

For the discussion above, it suffices to prove the necessity. Before proving,
we need the following lemma.

Denote $v_{i}$ the $(n-1)-$dimension vector corresponding to the orbit $ M_{i} $ with
$1$ at $j-$th entry for all $ j \in M_{i} $ and $0$ otherwise and $M$ the
set $\{v_{i} | 1 \leq i \leq r(n, K)\}. $ Let $ \Gamma=(\gamma_{st}) $ be an $(n-1)-$order
square matrix defined by $ \gamma_{st} = \zeta_{n}^{st}. $ Notice that $ \Gamma $
is invertible and $ \Gamma v_{i} \in F^{n-1}(1 \leq i \leq r(n, K)). $ Let $ V =
\{v \in \Q^{n-1} | \Gamma v \in F^{n-1}\} $ and $W$ be the vector space over $ \Q $ spanned
by the vectors belongs to $ M. $  We have the following results.
\par \vskip 0.4 cm

{\bf Lemma 1. } \quad  We have $ V = W $ and $ M $ is an orthogonal basis of $ V. $
\par  \vskip 0.2 cm

{\bf Proof. } \ Firstly, we claim that if $ v \in V, $ then $ \Gamma v \in W\bigotimes_{\Q}F. $
Let $ v = (w_{1}, w_{1}, \cdots, w_{n-1})^{T} \in V $ and
$ u = \Gamma v = (u_{1}, u_{1}, \cdots, u_{n-1})^{T} \in F^{n-1}. $ It is enough to show
$ u_{k} = u_{l} $ if $k, l$ in the same $ M_{i}. $ Recall that $ M_{i} $ is an orbit
obtained from $ G_{n}(p) $ for some $ p \in d(n) $ and $ n = pg. $ Let $ f(x) = u_{k} -
\sum _{h=1} ^{n-1}w_{h}x^{h} $  be a polynomial in $ F[x]. $ So $ f(\zeta_{n}^{k}) = 0 $
for $ u_{k} = \sum_{h=1}^{n-1}w_{h}\zeta_{n}^{kh}. $  Because $ k, l \in M_{i}, \zeta_{n}^{k}$
and $\zeta_{n}^{l}$ are two primitive $g-$th roots of unity. So $\zeta_{n}^{k}$
and $\zeta_{n}^{l} \in \Q(\zeta_{g}), $ which shows that
$ f(x) \in (F\bigcap \Q(\zeta_{g}))[x]. $ By the construction of $ M_{i}, $ there
exists an element $ \sigma \in  Gal(\Q(\zeta_{g})/F \bigcap \Q(\zeta_{g})) $ such
that $ \sigma(\zeta_{n}^{k}) = \zeta_{n}^{l}, $ which means that $\zeta_{n}^{k}$
is a conjugate element of $\zeta_{n}^{l}$ over $ F \bigcap \Q(\zeta_{g}) $ in $\Q(\zeta_{g}). $
Hence  $\zeta_{n}^{k}$ and $\zeta_{n}^{l}$ have the same minimal polynomial over
$ F \bigcap \Q(\zeta_{g}), $ which divides $ f(x). $ So $ f(\zeta_{n}^{l}) = 0, $
that is $ u_{l} = \sum_{h=1}^{n-1}w_{h}\zeta_{n}^{lh} = u_{k}. $
Using the result above, it is easy to get
$ \Gamma(V\bigotimes_{\Q}F) \subseteq W\bigotimes_{\Q}F, $
which shows that $ V\bigotimes_{\Q}F = W\bigotimes_{\Q}F. $
Therefore $ V = W $ and $ M $ is an orthogonal basis of $ V. $ The proof is completed.
  \quad \quad\quad $\Box$

Now it comes to prove the Theorem 1.
\par  \vskip 0.2 cm

{\bf Proof. } \ Consider the vector $ v \in \Q^{n-1} $ such that $1$ at
$j-$th entry for all $ j \in S $ and $0$ otherwise. Since $ D(n, S) $ is a circulant
digraph integral over $ K, \Gamma v = (\lambda_{1}, \lambda_{2},
\cdots, \lambda_{n-1})^{T} \in F^{n-1}. $ Hence $ v \in V $ and by Lemma1
$ v = \sum _{i=1} ^{r(n, K)}c_{i}v_{i} $ for some rational coefficients $c_{i}^{\prime}$s.
By the construction of $ M $ and $ v, S $ is a union of those $ M_{i}^{\prime} $s with
$ c_{i} = 1. $ The proof is completed. \quad \quad\quad \quad\quad $\Box$

By the Theorem 1, we obtain the following two corollaries.

Given a fixed pair $ (n, K), $ we have the partition $\{1, 2, \cdots, n-1\}
= \bigsqcup_{i=1}^{r(n, K)}M_{i}, $ where
$ r(n, K) = \sum_{p\in d(n)}[K \bigcap \Q(\zeta_{\frac{n}{p}}) : \Q]. $ Therefore there
are $ 2^{r(n, K)} $ distinct $ S $ for circulant digraph integral over $ K. $

{\bf Corollary 1. } \ There are at most $ 2^{r(n, K)} $ basic circulant digraphs
integral over $ K $ on $ n $ vertices.

Replacing $ K $ by the specific quadratic field $ \Q(i) $ and let $ n = 2^{t}k(t >2), $
it is easy to check that the partition of $\{1, 2, \cdots, n-1\} $ described in Part 3
of paper [XM] is the same as $ \bigsqcup_{i=1}^{r(n, K)}M_{i}. $ By Theorem 1,
we infer that the Conjecture 3.3 in [XM] is affirmative. Using the same notations of
[XM], we have the following corollary.

{\bf Corollary 2. } \ (Conjecture 3.3 in [XM]) Let $ D(\Z_{n}, S) $ be a
circulant digraph on $ n = 2^{t}k(t >2) $ vertices, then $ D(\Z_{n}, S) $ is Gauss integral
if and only if $ S $ is a union of the $M'$s.

\par  \vskip 1 cm
\hspace{-0.8cm} {\bf References }
\begin{description}

\item[[Ba]] K. Balinska, D.Cvetkovic, Z.Radosavljevic, S.Simic, D.Stevanovic,
A Survey on Integral Graphs, Univ. Beograd. Publ. Elektrotehn. Fak. Ser. Mat.
2002, 13 :42-65.

\item[[BC]] F. Bussemaker, D.Cvetkovic, There are exactly
connected, cubic, integral graphs, Univ. Beograd, Publ. Elektrotehn. Fak. Ser.
Mat. Fiz. 1976, 544: 43-48.

\item[[Bi]] N. Biggs, Algebraic Graph theory, North-Holland, Amsterdam, 1985.

\item[[EH]] F. Esser, F. Haraly, Digraphs with real and Gaussian spectra, \\
Discrete Appl. Math. 1980, 2: 113-124.

\item[[F]] K. Q. Feng, Algebraic Number Theory, (in Chinese), Science China Press, 2000.

\item[[GR]] C. Godsil, G. Royle, Algebraic Graph theory, New York: Springer-Verlag, 2001.

\item[[HS]] F. Haraly, A.Schwenk, Z.Radosavljevic, S.Simic, D.Stevanovic,
A Survey on  Which graphs have integral spectra, in:R. Bali, F. Haraly,(Eds.),
Graphs and Combinatorics, Berlin: Springer-Verlag, 1974, P.45.
2002, 13 :42-65.

\item[[L]] S. Lang, Algebra, Revised Third Edition, New
York: Springer-Verlag, 2002.

\item[[LL]] X. L. Li, G. N. Lin, On integral trees problems, KeXue TongBao,
1988, 33:802-806.

\item[[N]] J. Neukirch, Algebraic Number Theory, Berlin: Springer-Verlag, 1999.

\item[[S]] W. So, Integral circulant graphs,
Discrete Mathematics, 2005, 306: 153-158.

\item[[W]] L. C. Washington, Introduction to Cyclotomic Fields,
Second Edition, New York: Springer-Verlag, 2003.

\item[[XM]] Y. Xu, J. X. Meng, Gaussian integral circulant digraphs,
Discrete Mathematics, 2011, 311: 45-50.

\item[[Z]] H. R. Zhang, Constructing Integral Directed Graph by
Circulant Matrix Methods, J. Zhengzhou Univ. 2005, 37(4):28-34.

\end{description}

\end{document}